\documentclass[12pt]{article}
\newcommand{\be}{\begin{equation}}
      \newcommand{\ee}{\end{equation}}
      \newcommand{\ba}{\begin{eqnarray}}
       \newcommand{\ea}{\end{eqnarray}}
\newcommand{\ban}{\begin{eqnarray*}}
       \newcommand{\ean}{\end{eqnarray*}}

 \newcommand{\qed}{\hspace*{\fill}\rule{3mm}{3mm}\quad}

\newcommand{\sect}[1]{\section{#1} \setcounter{equation}{0}}

\begin{document}
 \newtheorem{defn}[lem]{Definition}
 \newtheorem{theo}[lem]{Theorem}
 \newtheorem{prop}[lem]{Proposition}
 \newtheorem{rk}[lem]{Remark}
 \newtheorem{ex}[lem]{Example}
 \newtheorem{note}[lem]{Note}
 \newtheorem{conj}[lem]{Conjecture}

\title{Curvature Estimates for the Ricci Flow II}
\author{Rugang Ye \\ {\small Department of Mathematics} \\
{\small University of California, Santa Barbara}\\}

\date{}
\maketitle

\sect{Introduction} 
\vspace{3mm}

In this paper we present several curvature estimates and convergence results for
solutions of the Ricci flow
\ba   \label{ricciflow}
\frac{\partial g}{\partial t}=-2Ric.
\ea
The said curvature estimates are
space-time analogues of the curvature estimates in [Ye3],
and depend on the smallness of
certain local space-time $L^{\frac{n+2}{2}}$ integrals of
the norm of the Riemann curvature tensor, where $n$ denotes the dimension of the manifold.
On the other hand, the said convergence results require finiteness
of space-time $L^{\frac{n+2}{2}}$ integrals of the norm of the Riemann curvature
tensor.  Note that these curvature estimates and convergence results
also serve as  characterizations of
blow-up singularities, see e.g. Remark 2 below.
 (The same can be said of the curvature estimates obtained in [Ye3].)

To formulate our results, we need some terminologies, most of which have already 
been 
used in [Ye3]. Consider a Riemannian manifold $(M, g)$ ($g$ denotes the metric) possibly
with boundary.   For convenience, we define the distance between two points of $M$ to be $\infty$, if they belong to 
two different connected components.  Consider a point $x\in M$.
If $x$ is in the interior of $M$, we define the distance $d(x, \partial M)=d_g(x, \partial M)$ to be
$\sup\{r>0: B(x,r)$ is compact and contained in the interior of  $ M\}$, where
$B(x, r)$ denotes the closed geodesic ball of center $x$ and radius $r$. If $M$ has
a boundary and $x \in \partial M$, then $d(x, \partial M)$ is the ordinary distance
from $x$ to $\partial M$ and equals zero. (For example, $d(x, \partial M)=\infty$ if
$M$ is closed.) 

For a family $U(t), 0\le t <T$ of open sets of $M$ for some $T>0$, we define its direct 
limit $\underline{\lim}_{t\rightarrow T}U(t)$ as follows. A point $x$ of $M$ 
lies in  $\underline{\lim}_{t\rightarrow T}U(t)$, if there is a neighborhood $U$ of $x$ 
and some $t\in [0, T)$ such that $U \subset U(t')$ for all $t' \in [t, T)$. \\

\noindent {\bf Notations} Let $g=g(t)$ be a family of metrics on $M$. Then $d(x, y, t)$ denotes the distance
between $x, y \in M$ with respect to the metric $g(t)$, and $B(x, r ,t)=B_g(x,r,t)$ denotes the closed geodesic ball of center $x\in M$ and radius $r$ with
respect to the metric $g(t)$. The volume of $B(x, r, t)$ with respect to $g(t)$ will 
often be denoted by $V(x,r,t)$ or $V_g(x, r, t)$. We shall often use $dq$ to denote $dvol_{g(t)}$.
These notations naturally extend when $g$ and (or) $t$ are replaced by 
other notations.

We set $\alpha_n=\frac{1}{40(n-1)}$, $\epsilon_0=\frac{1}{168}$ and
$\epsilon_1=\frac{\epsilon_0}{8\sqrt{1+2\alpha_n\epsilon_0^2}}$.
(These constants are not meant to be optimal. They can be improved
by closely examining the proofs.) \\

 We divide our results into several types.  In each type, the first theorem is
a local curvature estimate, the second theorem a global convergence result, and the
third theorem a local convergence result. In contrast to [Ye3], the results in this paper 
are also valid in dimension 2.  (Note that the Ricci flow is trivial in dimension 1.)  
\\
\hspace{1cm}

\noindent {\bf Type A} \\

Results of this type  involve straight (i.e. non-weighted) space-time $L^{\frac{n+2}{2}}$-integals
of the norm of the Riemann curvature tensor.  Theorem A-2 (the convergence result)
does not involve
any additional quantity or condition. Theorem A-1 and Theorem A-3
involve the condition of
$\kappa$-noncollapsedness, whose definition can be found in [P] and [Ye3].  Note that By [Theorem 4.1, P] and 
[Theorem A.1, Ye3], a smooth solution of the Ricci flow on $M \times [0, T)$
for a closed manifold $M$ and a finite $T$
is $\kappa$-noncollapsed on the scale $\rho$ for an arbitary positive number 
$\rho$, where $\kappa$ depends on
the initial metric and $T+\rho^2$.\\  \\

\noindent {\bf Theorem A-1} {\it  For each  positive number $\kappa$ and
each natural number $n\ge 2$
there are positive constants
$\delta_0=\delta_0(n, \kappa)$,
$C_0=C_0(n, \kappa)$ and $\sigma_0=\sigma_0(n, \kappa)$ depending only on
$n$ and $\kappa$ with the following property.
Let $g=g(t)$ be a smooth solution of the Ricci flow on
$M \times [0, T)$ for a manifold $M$ of dimension $n\ge 2$ and some
(finite or infinite) $T>0$, which is $\kappa$-noncollapsed  on
the scale $\rho$ for some $\kappa>0$ and $\rho>0$.
Consider $x_0 \in M$ and $0<r_0\le \rho$,which satisfy
$r_0<d_{g(t)}(x_0, \partial M)$ for
each $t \in [0, T)$. Assume that
\ba \label{smallRmA-1}
\int_0^T \int_{B(x_0, r_0, t)} |Rm|^{\frac{n+2}{2}}dqdt \le \delta_0.
\ea
Then we have
\ba \label{RmboundA-1}
|Rm|(x, t) \le \alpha_n t^{-1}+(\epsilon_0 r_0)^{-2}
\ea
whenever $t \in (0, T)$ and $d(x_0, x, t)<\epsilon_0 r_0$, and 
\ba \label{RmboundintA-1}
|Rm|(x, t) \le C_0 \max\{\frac{1}{r_0^{2}}, \frac{1}{t}\}(\int_0^t \int_{B(x, \frac{1}{2}r(t), s)} |Rm|^{\frac{n+2}{2}}dqds)^{\frac{2}{n+2}},
\ea
whenever $0<t<T$ and $d(x_0, x, t)\le \frac{1}{2}r(t)$, where $r(t)=\epsilon_1 \min\{r_0, \sqrt{t}\}.$ 
(Obviously, the estimates (\ref{RmboundA-1}) and (\ref{RmboundintA-1})
hold on $[0, T]$ provided that $T$ is finite and the assumptions hold on
$[0, T]$. This remark also applies to the results below.)} \\

Note that the constant $\delta_0$ depends on $n$ decreasingly and depends on $\kappa$ increasingly, i.e. 
$\delta_0(n, \kappa)$ is a decreasing function of $n$ and an increasing function of $\kappa$.  In contrast, the 
constant $C_0$ depends on $n$ increasingly and depends on $\kappa$ decreasingly.  
The dependences of the constants $\delta_0$ and $C_0$ in Theorem B-1 and Theorem C-1 are of similar nature. \\

\noindent {\bf Remark 1} Theorem A-1 is optimal in the sense that
if we replace $\frac{n+2}{2}$ by a smaller exponent, then
the conclusion fails to hold. This is demonstrated by the example of
the evolving sphere. This remark also applies to the results below. \\

\noindent {\bf Theorem A-2} {\it Let $g=g(t)$ be a smooth solution of
the Ricci flow on $M \times [0, T)$ for an n-dimensional
closed manifold $M$ of dimension $n \ge 2$ and some finite $T>0$.
Assume
\ba \label{finiteRmA-2}
\int_0^T \int_M |Rm|^{\frac{n+2}{2}} dq dt <\infty.
\ea
Then $g(t)$ converges smoothly to a smooth metric on $M$ as
$t\rightarrow T$.
Consequently, $g(t)$ extends to a smooth solution
of the Ricci flow over $[0,T']$ for some $T'>T$.
} \\

\noindent {\bf Remark 2} Theorem A-2 can be rephrased as follows:
A (global or local) solution of the Ricci flow  blows up at $T$,
if and only if the space-time integral of $|Rm|^{\frac{n+2}{2}}$ up to $T$ is infinite. This
can be used to analyse blow-ups of the Ricci flow. For example, careful rescalings
produce blow-up limits with the special feature of infinite space-time integral of $|Rm|^{\frac{n+2}{2}}$.  This remark also applies to the results below. \\

\noindent {\bf Theorem A-3} {\it
Let $g=g(t)$ be a smooth solution of the Ricci flow on
$M \times [0, T)$ for a manifold $M$ of dimension $n\ge 2$ and some
finite  $T>0$, which is $\kappa$-noncollapsed  on
the scale of $\rho$ for some $\kappa>0$ and $\rho>0$.
Consider $x_0 \in M$ and $0<r_0\le \rho$, which satisfy
$r_0\le diam_{g(t)}(M)$ and $r_0<d_{g(t)}(x_0, \partial M)$ for
each $t \in [0, T)$. Assume that
\ba \label{finiteRmA-3}
\int_0^T \int_{B(x_0, r_0, t)} |Rm|^{\frac{n+2}{2}}dq dt<\infty.
\ea
If $T<\infty$, then $g(t)$ converges smoothly to a smooth metric $g(T)$ on
the direct limit $\underline{\lim}_{t \rightarrow T} \mathring{B}(x_0, \epsilon_0,$
$t)$. Moreover, $\mathring{B}(x_0, \epsilon_0 r_0, T)=\underline{\lim}_{t \rightarrow T} \mathring{B}(x_0, \epsilon_0r_0,
t).$  
}\\

\hspace{1cm}

\noindent {\bf Type B} \\

Results of this type do not involve the condition of $\kappa$-noncollapsedness.
Instead, they employ space-time $L^{\frac{n+2}{2}}$
integrals of the norm of $Rm$
over balls of varying center and radius measured against a volume ratio. \\

\noindent {\bf Theorem B-1} {\it
For each natural number $n\ge 2$ there are  positive constants $
\delta_0=\delta_0(n)$, $C_0=C_0(n)$ and $\sigma_0=\sigma_0(n)$ depending only on $n$ with the following property.
Let $g=g(t)$ be a smooth solution of the Ricci flow on
$M \times [0, T)$ for a connected manifold $M$ of dimension $n \ge 2$ and some
(finite or infinite) $T>0$.
Consider $x_0 \in M$ and $r_0>0$,which satisfy
$r_0\le diam_{g(t)}(M)$ and $r_0<d_{g(t)}(x_0, \partial M)$ for
each $t \in [0, T)$.  Assume that
\ba \label{smallRmB-1}
\int_0^T \sup \limits_{x \in
B(x_0, \frac{r_0}{2}, t)}\frac{r^n}{V(x, r, t)} \int\limits_{B(x, r, t)} |Rm|^{\frac{n+2}{2}}(\cdot, t)dq(\cdot, t)dt \le \delta_0,
\ea
for all $0<r\le \frac{r_0}{2}$.  Then we have
\ba \label{RmboundB-1}
|Rm|(x, t) \le \alpha_n t^{-1}+(\epsilon_0 r_0)^{-2}
\ea
whenever $t \in (0, T)$ and $d(x_0, x, t)<\epsilon_0 r_0$, and 
\ba \label{RmboundintB-1}
|Rm|(x, t) \le \frac{C_0}{ r(t)^{2}}
\left(\int_0^t \frac{r(t)^n}{V(x,\frac{1}{9}r(t), s)}\int\limits_{B(x, \frac{1}{9}r(t), s)} |Rm|^{\frac{n+2}{2}}(\cdot,s)dq(\cdot, s)
ds\right)^{\frac{2}{n+2}},
\ea
whenever $0<t<T$ and $d(x_0, x, t)\le \frac{1}{2}r(t)$, where $r(t)= \epsilon_1 \min\{r_0, \sqrt{t}\}$.
} \\

\noindent {\bf Theorem B-2} {\it Let $g=g(t)$ be a smooth solution
of the Ricci flow  on $M \times [0, T)$ for a  closed
manifold $M$ of dimension $n \ge 2$ and some finite $T>0$. Assume 
\ba \label{finiteRmB-2} 
\int_0^T \sup\limits_{x \in M, 0<r \le d_{\mu}(r_0, g(t))} 
\left(\frac{r^n}{V(x, r, t)} \int_{B(x, r, t)}
|Rm|^{\frac{n+2}{2}}(\cdot, t)dq(\cdot, t) \right) dt
 <\infty
\ea 
for some $r_0>0$ and $0<\mu \le 1$, where $d_{\mu}(r_0, g(t)) =
\min\{r_0, \mu \cdot diam_{g(t)}(M)\}$.   Then $g(t)$ converges smoothly
to a smooth metric on $M$ as $t \rightarrow T$. Consequently,
$g(t)$ extends to a smooth solution of the Ricci flow over
$[0,T']$ for some $T'>T$.
} \\

\noindent {\bf Remark 3} We obtain two interesting
special cases of this theorem when we replace the condition
$0<r<\min\{r_0, diam_{g(t)}(M)\}$  in (\ref{finiteRmB-2}) by $0<r<r_0$
or $0<r<\mu diam_{g(t)}(M)$. \\

\noindent {\bf Theorem B-3} {\it
Let $g=g(t)$ be a smooth solution of the Ricci flow on
$M \times [0, T)$ for a manifold $M$ of dimension $n \ge 2$ and some
finite $T>0$.
Consider $x_0 \in M$ and $r_0>0$, which satisfy
$r_0\le diam_{g(t)}(M)$ and $r_0<d_{g(t)}(x_0, \partial M)$ for
each $t \in [0, T)$.
Assume
\ba \label{finiteRmB-4}
\sup\limits_{0<r \le \frac{r_0}{2}}\int_0^T \left( \sup\limits_{x\in B(x_0, \frac{r}{2}, t)}  \frac{r^n}{V(x, r, t)} \int_{B(x, r, t)} |Rm|^{\frac{n+2}{2}}(\cdot,t)dq(\cdot, t) \right)dt<\infty.
\ea
Then $g(t)$ converges smoothly to a smooth metric $g(T)$ on
the direct limit $\underline{\lim}_{t \rightarrow T} \mathring{B}(x_0,$
 $\epsilon_0 r_0,t)$. Moreover, $\mathring{B}(x_0, \epsilon_0 r_0, T)=\underline{\lim}_{t \rightarrow T} \mathring{B}(x_0, \epsilon_0,
t).$
}\\

\hspace{1cm}

\noindent {\bf Type C} \\

Results of this type do not involve the condition of $\kappa$-noncollapsedness,
and employ only the space-time
$L^{\frac{n+2}{2}}$ integrals of the
norm of the Riemann curvature tensor over balls of a {\it fixed center} and 
a {\it fixed radius}, measured against a volume ratio. But a lower bound for the
Ricci curvature is assumed. \\

\noindent {\bf Theorem C-1} {\it
For each natural number $n\ge 2$ there are  positive constants $
\delta_0=\delta_0(n)$ and $C_0=C_0(n)$ depending only on $n$ with the following property.
Let $g=g(t)$ be a smooth solution of the Ricci flow on
$M \times [0, T)$ for a  connected manifold $M$ of dimension $n \ge 2$ and some
finite $T>0$.
 Consider $x_0 \in M$ and $r>0$, which satisfy
$r_0\le diam_{g(t)}(M)$ and $r_0<d_{g(t)}(x_0, \partial M)$ for
each $t \in [0, T)$. Assume that
\ba \label{ricciboundC-1}
Ric(x, t) \ge -\frac{n-1}{r_0^2}g(x,t)
\ea
whenever $t \in [0, T)$ and $d(x_0, x, t)\le r_0$ ($g(x, t)=g(t)(x)$ and 
$Ric(x, t)$ is the Ricci tensor of $g(t)$ at $x$), and that
\ba \label{smallRmC-1}
\int_0^T \frac{r_0^n}{V(x_0, r_0, t)} \int_{B(x_0, r_0, t)} |Rm|^{\frac{n+2}{2}}(\cdot, t)dq(\cdot, t) dt \le \delta_0.
\ea
Then we have
\ba \label{RmboundC-1}
|Rm|(x, t) \le \alpha_n t^{-1}+(\frac{1}{2}\epsilon_0 r_0)^{-2}
\ea
whenever $t \in (0, T)$ and $d(x_0, x, t)<\frac{1}{2}\epsilon_0 r_0$,
and
\ba \label{RmboundintC-1}
|Rm|(x, t) \le \bar C_0 r(t)^{-2}
\left(\int_0^t \frac{r_0^n}{V(x_0, r_0, s)}\int\limits_{B(x_0, r_0, s)} |Rm|^{\frac{n+2}{2}}(\cdot,s)dq(\cdot, s)
ds\right)^{\frac{2}{n+2}}
\ea
whenever $0<t<T$ and $d(x_0, x, t)\le \frac{1}{2}r(t)$, where $r(t)=\epsilon_1 \min\{\frac{1}{2}r_0, \sqrt{t}\}$.
We also have
\ba \label{RmboundintC-2}
|Rm|(x, t) \le C_0 r(t)^{-2}
\left(\int_0^t \frac{r(t)^n}{V(x, \frac{1}{9}r(t), s)}\int\limits_{B(x, \frac{1}{9}r(t), s)} |Rm|^{\frac{n+2}{2}}(\cdot,s)dq(\cdot, s)
ds\right)^{\frac{2}{n+2}}
\ea
whenever $0<t<T$ and $d(x_0, x, t)\le \frac{1}{2}r(t)$. 
} \\

\noindent {\bf Theorem C-2} {\it Let $g=g(t)$ be a smooth solution of
the Ricci flow on $M \times [0, T)$ for a closed
manifold $M$ of dimension $n \ge 2$ and some finite $T>0$. Assume that (\ref{ricciboundC-1})
holds for all $x\in M$ and $t \in [0,T)$, and
\ba \label{finiteRmC-2}
\sup\limits_{x\in M} \int_0^T \frac{1}{V(x, r_0, t)} \int_{B(x, r_0, t)} |Rm|^{\frac{n+2}{2}}(\cdot,t) dq(\cdot, t) dt <\infty
\ea
for some $r_0>0$. 
Then $g(t)$ converges smoothly to a smooth metric on $M$ as
$t\rightarrow T$. Consequently, $g(t)$ extends to a smooth solution
of the Ricci flow over $[0,T']$ for some $T'>T$.
 } \\

\noindent {\bf Theorem C-3} {\it
For each natural number $n\ge 2$ there is a positive constant $
\delta_0=\delta_0(n)$ depending only on $n$ with the following property.
Let $g=g(t)$ be a smooth solution of the Ricci flow on
$M \times [0, T)$ for a manifold $M$ of dimension $n \ge 2$ and some
finite $T>0$.
Consider $x_0 \in M$ and $r_0>0$,
which satisfy
$r_0\le diam_{g(t)}(M)$ and $r_0\le d_{g(t)}(x_0, \partial M)$ for
each $t \in [0, T)$.
Assume that (\ref{ricciboundC-1}) holds whenever $0\le t <T$ and
$d(x_0, x, t)<r_0$, and that
\ba \label{finiteRmC-3}
\int_0^T \frac{1}{V(x_0, r_0, t)} \int_{B(x_0, r_0, t)} |Rm|^{\frac{n+2}{2}}(\cdot,t)dq(\cdot, t) dt <\infty.
\ea
Then $g(t)$ converges smoothly to a smooth metric $g(T)$ on
the direct limit $\underline{\lim}_{t \rightarrow T} \mathring{B}(x_0,$
 $\epsilon_0 r_0,t)$. Moreover, $\mathring{B}(x_0, \epsilon_0 r_0, T)=\underline{\lim}_{t \rightarrow T} \mathring{B}(x_0, \epsilon_0,
t).$
}\\

\vspace{1cm}

Note that the condition $r_0 \le diam_{g(t)}(M)$ appears in Theorem B-1 and Theorem C-1, but not 
in Theorem A-1.  \\

Now we discuss extensions of the above results. First, Theorem A-2, Theorem B-2 and Theorem C-2  extend to noncompact manifolds  under an additional assumption of $\kappa$-noncollapsedness. 
Suitable extensions of Theorem A-2, Theorem A-3, Theorem B-2, Theorem B-3, Theorem C-2, and Theorem C-3 
also hold true 
in the case $T=\infty$.

Theorem A-1, Theorem B-1 and Theorem C-1 extend  to 
the modified Ricci flow 
\ba \label{chi}
\frac{\partial g}{\partial t}=-2Ric+\lambda(g, t) g
\ea
with a scalar function $\lambda(g, t)$ independent of $x \in M$. 
(The volume-normalized 
Ricci flow 
\ba \label{volumenormalize}
\frac{\partial g}{\partial t}=-2Ric+\frac{2}{n} {\hat R} g
\ea
on a closed manifold, with $\hat R$ denoting the average scalar curvature, is  an example of the modified Ricci flow.)
We present two extensions which are analogous to Extension I and Extension II in [Ye3].\\

\noindent {\bf Extention I} {\it Theorem A-1, Theorem B-1 and Theorem C-1 hold true for the modified Ricci flow (\ref{chi}), with the modification that 
the constants $\delta_0$ and $C_0$ in each theorem 
depend in addition on $r_0^2 |\min\{\inf_{[0, T)} \lambda, 0\}|$ which is assumed to be finite.
(In other words, $\delta_0$ and $C_0$ depend in addition on a nonpositive lower bound of $r_0^2 \lambda$. )
 } \\

\noindent {\bf Extention II} {\it Theorem A-1, Theorem B-1 and Theorem C-1 in the case $T<\infty$ hold true for the 
modified Ricci flow (\ref{chi}), with the modification that the constants $\delta_0$ and $C_0$ in 
each theorem depend in addition on $|\min\{\inf_{0\le t_1 <t_2 <T} \int_{t_1}^{t_2} \lambda, 0\}|$ 
which is assumed to be finite.  (In other words, $\delta_0$ and $C_0$ depend in addition on a nonpositive 
lower bound of $\int_{t_1}^{t_2} \lambda$.) } \\ 

In both extentions, the dependence of $\delta_0$ is decreasing, and the dependence of $C_0$ is increasing.   Extention I can be proved by directly adapting the proofs of Theorem A-1, Theorem B-1 and Theorem C-1. Extention II 
can be proved by converting the modified Ricci flow into the Ricci flow, applying Theorem A-1, Theorem B-1 and 
Theorem C-1, and then converting the obtained estimates back to the modified Ricci flow.  
(Such an argument 
can be found in the proof of Theorem B-2.)

Suitable extensions of Theorem A-2, Theorem A-3, Theorem B-2, Theorem B-3, Theorem C-2, and Theorem C-3
also hold true for the modified Ricci flow. 

We would also like to point out that the results in this paper can be extended to many other evolution  equations in various ways. 

The curvature estimates in this paper were obtained some time ago.\\

\sect{Curvature Estimates}

\vspace{3mm}

\subsection{Type A} 
\vspace{3mm}

In this subsection we present the proof of Theoem A-1, which 
is divided into two parts. \\

\noindent {\bf Proof of the estimate (\ref{RmboundA-1})} \\

The proof is similar to [Proof of the estimate (1.3), Ye3]. To make the proof
clear, we'll repeat some arguments  in [Ye3]. By rescaling, we can assume $r_0=1$.
Assume that the estimate (\ref{RmboundA-1}) does not hold. Then we can find for each
$\epsilon>0$
a Ricci flow solution $g=g(t)$ on $M\times [0, T)$ for some $M$ and $T>0$
with the properties as postulated in the statement of the theorem, such that $|Rm|(x, t)>\alpha_n t^{-1}+\epsilon^{-2}$ for some
$(x, t) \in M \times (0, T)$ satisfying $d(x_0, x, t) <\epsilon$.

We denote by $M_{\alpha_n}$ the set of pairs $(x, t)$ such that
$|Rm|(x, t)\ge \alpha_n t^{-1}$.
For an arbitary  positive number $A>1$ such that
$(2A+1)\epsilon \le \frac{1}{2}$,
we choose as in [Proof of Theorem 10.1, P] and [Proof of Theorem A, Ye3] a point $(\bar x, \bar t) \in M_{\alpha_n}$ with
$0<\bar t \le \epsilon^2, d(x_0, \bar x, \bar t)<(2A+1)\epsilon$,
such that $|Rm|(\bar x, \bar t) > \alpha_n \bar t^{-1}+\epsilon^{-2}$
and
\ba \label{4RmA}
|Rm|(x, t)\le 4|Rm|(\bar x, \bar t)
\ea
whenever
\ba \label{4Rmcondition}
(x, t) \in M_{\alpha_n}, 0 <t \le \bar t, d(x_0, x, t)\le d(x_0, \bar x, \bar t)
+A|Rm|(\bar x, \bar t)^{-\frac{1}{2}}.
\ea

We set $Q=|Rm|(\bar x, \bar t)$. By [Proof of Theorem A, Ye3], the following
two claims hold. \\

\noindent {\bf Claim 1}  If
\ba \label{new4Rmcondition}
\bar t -\frac{1}{2}\alpha_n Q^{-1} \le t \le \bar t,
d(\bar x, x, \bar t) \le \frac{1}{10} AQ^{-\frac{1}{2}},
\ea
then
\ba \label{distance1}
d(x_0,x,t) \le d(x_0,\bar x,\bar t)+\frac{1}{2} AQ^{-\frac{1}{2}}.
\ea
\\

\noindent {\bf Claim 2} If $(x, t)$ satisfies (\ref{new4Rmcondition}), then
the estimate (\ref{4RmA}) holds. \\

An implication of  Claim 1 is
\ba \label{distance2}
d(x_0, x, t) \le (\frac{5}{2}A+1)\epsilon
\ea
for $(x, t)$ satisfying (\ref{new4Rmcondition}).

Now we take $\epsilon=\frac{1}{42}$ and $A=10$ as in [Ye3].  Then $\frac{1}{10}A<1$ and
$(\frac{5}{2}A+1)\epsilon=1$. So (\ref{distance2}) implies
\ba
B(\bar x, Q^{-\frac{1}{2}}, \bar t) \subset B(x_0, 1, t)
\ea
for $t \in [\bar t-\frac{1}{2}\alpha_n Q^{-1}, \bar t]$, and hence
the condition (\ref{smallRmA-1}) leads to
\ba 
\int_{\bar t-\frac{1}{2}\alpha_n Q^{-1}}^{\bar t}\int_{B(\bar x, Q^{-\frac{1}{2}}, \bar t)} |Rm|^{\frac{n+2}{2}}(\cdot, t) dq(\cdot, t)dt \le \delta_0.
\ea
   Moreover, Claim 2 implies
that the estimate (\ref{4RmA}) holds on $B(\bar x, Q^{-\frac{1}{2}}, \bar t)
\times [\bar t-\frac{1}{2}\alpha_n Q^{-1}, \bar t]$.
 As in [Ye3], we shift $\bar t$ to the time origin
and rescale $g$ by the factor $Q$ to obtain a Ricci flow solution $\bar g(t)=Qg(\bar t+Q^{-1}t)$
on $M \times [-\frac{1}{2}\alpha_n, 0]$.
Then we have for $\bar g$
\ba \label{point1}
|Rm|(\bar x, 0)=1,
\ea
\ba \label{new4Rm1}
|Rm|(x, t) \le 4
\ea
whenever
$$
-\frac{1}{2}\alpha_n \le t \le 0,
d(\bar x, x, 0) \le 1,
$$
and
\ba \label{proofA-1}
\int_{-\frac{1}{2}\alpha_n}^0\int_{B(\bar x, 1, 0)}|Rm|^{\frac{n+2}{2}}(\cdot, t)dq(\cdot ,t)dt \le \delta_0.
\ea
Moreover, $\bar g$ is $\kappa$-noncollpased on the scale 
$Q^{\frac{1}{2}} \rho \ge 1$. We also have 
\ba \label{newdistance}
d_{\bar g(0)}(\bar x, \partial M) > \frac{1}{2}.
\ea

As in [Ye3], we apply the above properties to deduce $C_{S, 2, \bar g(0)}(B(\bar x, \rho(n, \kappa), 0)) 
\le C_1(n)$ and then
\ba
C_{S,2, \bar g (t)}(B(\bar x, \rho(n, \kappa), 0)))\le C_2(n)
\ea
for $t \in [-\frac{1}{2}\alpha_n, 0]$, where $C_1(n)>0$ and $ C_2(n)>0$ 
depend only on $n$, and 
$0<\rho(n, \kappa)\le \frac{1}{16}$ depends only on $n$ and $\kappa$.

As in [Ye3] we have  $B(\bar x, \rho(n, \kappa), 0)
\subset B(\bar x, 1, t)$ for $t \in [-\bar \alpha_n, 0]$, where
$\bar \alpha_n\le \frac{1}{2}\alpha_n$ is a positive constant depending only on $n$.
It follows that
\ba \label{smallA}
\int_{-\bar \alpha_n}^0\int_{B(\bar x, \rho(n, \kappa), 0)} |Rm|^{\frac{n+2}{2}}(\cdot, t)dq(\cdot, t)dt \le
\delta_0.
\ea

As in [Ye3], we now appeal to the differential inequality
\ba \label{curvaturenorm}
\frac{\partial}{\partial t}|Rm| \le \Delta |Rm| + c(n) |Rm|^2
\ea
for a positive constant $c(n)$ depending only on $n$.
 On account of (\ref{new4Rm1}), (\ref{sobolevA}) and
(\ref{smallA}) we can apply [Theorem 2.1, Ye3] to (\ref{curvaturenorm}) with $p_0=\frac{n+2}{2}$  to deduce
\ba \label{RmfinalA-1}
|Rm|(\bar x, 0)
&\leq& (1+\frac{2}{n})^{\frac{2\sigma_n}{n+2}} C_2(n)^{\frac{2n}{n+2}} C_3(n, \kappa)
\Bigl(\int^0_{-\bar \alpha_n}\int\limits_{B(\bar x, \rho(n, \kappa),0)} |Rm|^{\frac{n+2}{2}}(\cdot, t)
dq(\cdot, t)dt
\Bigr)^{\frac{2}{n+2}} \nonumber \\
&\le& (1+\frac{2}{n})^{\frac{2\sigma_n}{n+2}} C_2(n)^{\frac{2n}{n+2}} C_3(n, \kappa)
\delta_0
^{\frac{2}{n+2}},
\ea
where
\ban
C_3(n, \kappa)=2c(n)(n+2)+4n(n-1) +\frac{n(n+2)^2}{8\bar \alpha_n}
+
\frac{(n+2)^2}{4\rho(n, \kappa)^2}e^{8(n-1) \bar \alpha_n}.
\ean
We deduce $|Rm|(\bar x, 0) \le \frac{1}{2}$, provided that we define
\ban
\delta_0 &=&2^{-\frac{n+2}{2}}(1+\frac{2}{n})^{-\sigma_n} C_2(n)^{-n} C_3(n, \kappa)^{-\frac{n+2}{2}}.
\ean
But this contradicts (\ref{point1}). Hence the estimate (\ref{RmboundA-1})
 has been proven.
\\

\noindent {\bf Proof of the estimate (\ref{RmboundintA-1})} \\

The proof of [(1.4), Ye3] in [Ye3] carries over.  Consider a fixed $t_0 \in (0, T)$. 
We set $\bar g(t)=\lambda_n t_0^{-1} g(\lambda_n^{-1} t_0 t$ if $r_0^2 \ge t_0$ and 
$\bar g(t) =\lambda_n r_0^{-2} g(\lambda_n^{-1} r_0^2 t)$ if $r_0^2 \le t_0$.  We replace 
the exponent $\frac{n}{2}$ in the said proof (namely in the first line of [(3.27), Ye3]) by 
$\frac{n+2}{2}$ and 
deduce for $\bar g$
\ba
|Rm|(x, \lambda_n) \le C_0(n, \kappa)
\left(\int_{\lambda_n-\bar \alpha_n}^{\lambda_n} \int_{B(x, \rho(n, \kappa), \lambda_n)}|Rm|^{\frac{n+2}{2}}(\cdot, t)dq(\cdot, t)dt\right)^{\frac{2}{n+2}}
\ea
for all $x\in B(x_0, \frac{1}{2}, \lambda_n)$, where $\lambda_n>0$ and $0<\bar \alpha_n<\frac{\lambda_n}{2}$ depend only on $n$, and $C_0(n, \kappa)>0$ and $\rho(n, \kappa)>0$ depend only on 
$n$ and $\kappa$.  Scaling back to $g$ we then arrive at the desired estimate (\ref{RmboundintA-1}). \\
\qed 
\\

\subsection{Type B} 

\vspace{3mm} 

In this subsection we present the proof of Theorem B-1, which is also 
divided into two parts. \\

\noindent {\bf Proof of the estimate (\ref{RmboundB-1})} \\

This is similar to the proof of [(1.6), Ye3] in [Ye3]. We'll repeat 
most arguments there for clarify.
Assume that the estimate (\ref{RmboundB-1}) fails to hold.
Then we carry out the same construction as in the proof of Theorem A-1.
Again we assume $r_0=1$ and choose  $\epsilon=\frac{1}{84}$ and $A=10$. 
We deal with the rescaled flow $\bar g$ and all quantities will be associated with 
$\bar g$.  By (\ref{smallRmB-1}) we have for $\bar g$, in place of (\ref{proofA-1})
\ba \label{proofB-1}
\int_{-\frac{1}{2}\alpha_n}^0 \frac{1}{V(\bar x, r, t)} \int_{B(\bar x, r, t)} |Rm|^{\frac{n+2}{2}}
dqdt  \le  \frac{\delta_0}{r^n}
\ea
whenever $0<r\le \frac{1}{2}$. Moreover, we have 
\ba
diam_{\bar g(t)} (M) \ge 1
\ea
for all $t \in [-\frac{1}{2}\alpha_n, 0]$ and 
\ba
d_{\bar g(t)}(\bar x, \partial M) \ge \frac{1}{3}
\ea
for all $t \in [-\alpha_n^*, 0]$, where $0<\alpha_n^* \le \frac{1}{2}\alpha_n$ 
depends only on $n$.  
 As before, we also have for $\bar g$
\ba \label{pointB-1}
|Rm|(\bar x, 0)=1
\ea
and
\ba \label{4RmB-1}
|Rm|(x, t) \le 4
\ea
whenever
\ba
-\frac{1}{2}\alpha_n\le t \le 0, d(\bar x, x, 0) \le 1.
\ea
Moreover, we have
\ba
diam_{\bar g(t)}(M) \ge 1
\ea 
for all $t \in [-\frac{1}{2}\alpha_n, 0]$ and 
\ba
d_{\bar g(t)}(\bar x, \partial M) \ge \frac{1}{3}
\ea
for all $t \in [-\alpha_n^*, 0]$, where $0<\alpha_n^*
\le \frac{1}{2}\alpha_n$ depends only on $n$. 

As in [Ye3] we derive from the above properties 
\ba \label{sobolevB-1}
C_{S, 2, \bar g(t)}(B(\bar x, \frac{1}{12}, t)) \le \frac{C_1(n)}{V(\bar x, \frac{1}{12}, t)^{\frac{1}{n}}}
\ea
for $t \in [-\alpha_n^*, 0]$, with a positive constant $C_1(n)$
depending only on $n$.  Moreover, we have
\ba \label{includeB-1}
B(\bar x, \frac{1}{14}, t_1) \subset  B(\bar x, \frac{1}{13}, t_2) \subset
B(\bar x, \frac{1}{12}, t_3)
\ea
for all $t_1, t_2$ and $t_3 \in [-\bar \alpha_n, 0]$, with a positive constant $\bar \alpha_n
\le \alpha_n^*$ depending only on $n$.
  Consequently, we derive from (\ref{proofB-1}) and (\ref{sobolevB-1})
\ba \label{smallintegralB-1-2}
\int_{-\bar \alpha_n}^0 \frac{1}{V(\bar x, \frac{1}{13}, t))}\int_{B(\bar x, \frac{1}{14}, 0)} |Rm|^{\frac{n+2}{2}}(\cdot, t)dq(\cdot, t)dt  \le (13)^n \delta_0
\ea
and
\ba \label{sobolevB-1-2}
C_{S, 2, \bar g(t)}(B(\bar x, \frac{1}{14}, 0)) \le \frac{C_1(n)}{V(\bar x, \frac{1}{12}, t))^{\frac{1}{n}}}
\ea
for all $t\in [-\bar \alpha_n, 0]$. Moreover, (\ref{includeB-1}) combined
with (\ref{4RmB-1}) leads via the Ricci flow equation to
\ba \label{volumeB-1-2}
\min\limits_{-\bar \alpha_n \le t \le 0} V(\bar x, \frac{1}{12}, t)
\ge e^{-4n(n-1)\bar \alpha_n} \max\limits_{-\bar \alpha_n\le t \le 0} V(\bar x, \frac{1}{13}, t)
\ea
for each $t \in [-\bar \alpha_n, 0]$.
Now we apply [Theorem 2.1,Ye3] to deduce
\ba \label{RmfinalB-1}
|Rm|(\bar x, 0)
&\leq& \frac{(1+\frac{2}{n})^{\frac{2\sigma_n}{n+2}} C_1(n)^{\frac{2n}{n+2}}C_2(n)}{\min\limits_{-\bar \alpha_n \le t \le 0} V(\bar x, \frac{1}{12}, t)^{\frac{2}{n+2}}}
\Bigl(\int^0_{-\bar \alpha_n}\int\limits_{B(\bar x, \frac{1}{14},0)} |Rm|^{\frac{n+2}{2}}
\Bigr)^{\frac{2}{n+2}} \nonumber \\
&\le& \frac{(1+\frac{2}{n})^{\frac{2\sigma_n}{n+2}}
C_1(n)^{\frac{2n}{n+2}}C_2(n)}{\min\limits_{-\bar \alpha_n \le t \le 0} V(\bar x, \frac{1}{12}, t)^{\frac{2}{n+2}}}  \max\limits_{-\bar \alpha_n
\le t \le 0} V(\bar x, \frac{1}{13}, t)^{\frac{2}{n+2}} \nonumber
\\ && \cdot
(\int_{-\bar \alpha_n}^0\frac{1}{V(\bar x, \frac{1}{13}, t)} \int\limits_{B(\bar x, \frac{1}{14},0)}
 |Rm|^{\frac{n+2}{2}}
\Bigr)^{\frac{2}{n+2}} \nonumber \\
 \nonumber \\
&\le& (13)^{\frac{2n}{n+2}}(1+\frac{2}{n})^{\frac{2\sigma_n}{n+2}} C_1(n)^{\frac{2n}{n+2}} C_2(n)
\delta_0
^{\frac{2}{n+2}}e^{\frac{2n(n-1)}{n+2} \bar\alpha_n},
\ea
with a suitable positive constant $C_2(n)$ depending only on $n$. Choosing
\ban
\delta_0=\frac{1}{(13)^n}(1+\frac{2}{n})^{-\sigma_n} C_1(n)^{-n} C_2(n)^{-\frac{n+2}{2}}
 e^{-n(n-1)\bar \alpha_n}
\ean
we then obtain $|Rm|(\bar x, 0)\le \frac{1}{2}$, contradicting
(\ref{pointB-1}). \\  \qed
\\

\noindent {\bf Proof of the estimate (\ref{RmboundintB-1})} \\

The arguments in the proof of [(1.7), Ye3] in [Ye3] carry over.  Consider a fixed $t_0 \in (0, T)$. 
We set $\bar g(t)=\lambda_n t_0^{-1} g(\lambda_n^{-1} t_0 t)$ if $r_0^2 \ge t_0$ and 
$\bar g(t) =\lambda_n r_0^{-2} g(\lambda_n^{-1} r_0^2 t)$ if $r_0^2 \le t_0$. 
We employ radii $\frac{1}{8}, \frac{1}{9}$ and $\frac{1}{10}$ as in the said proof, which play the role of the above radii $\frac{1}{14}, \frac{1}{13}$ and 
$\frac{1}{12}$. We deduce for $\bar g$ as in  (\ref{RmfinalB-1})  
\ba
|Rm|(x, \lambda_n) &\le& C_3(n)  \left( \int_{\lambda_n-\bar \alpha_n}^{\lambda_n} \frac{1}{V(x, \frac{1}{9}, t)} \int_{B(x, \frac{1}{10}, \lambda_n)} |Rm|^{\frac{n+2}{2}}(\cdot, t)dq(\cdot, t)\right)
^{\frac{2}{n+2}} \nonumber \\
&\le& C_0(n)  \left( \int_{\lambda_n-\bar \alpha_n}^{\lambda_n} \frac{(\frac{1}{9})^n}{V(x, \frac{1}{9}, t)} \int_{B(x, \frac{1}{9}, t)} |Rm|^{\frac{n+2}{2}}(\cdot, t)dq(\cdot, t)\right)
^{\frac{2}{n+2}},
\ea
whenever $d(x_0, x, \lambda_n)
\le \frac{1}{2}$, where $C_3(n)$ and $\bar \alpha_n$ and $C_0(n)$ depend only on $n$, and $C_0(n)=9^nC_3(n)$.   Scaling back to $g$ we then arive at the desired 
estimate (\ref{RmboundintB-1}). 
\\ \qed \\

\subsection{Type C}

\vspace{3mm}

\noindent {\bf Proof of Theorem C-1} \\

We establish the condition (\ref{smallRmB-1}). Then the estimate 
(\ref{RmboundC-1}) 
follows from Theorem B-1. This is similar to the proof of [Theorem C, Ye3] 
in [Ye3]. By rescaling we can assume $r_0=1$. Then
(\ref{ricciboundC-1}) becomes \ba \label{ricciboundC-1new} Ric \ge
-(n-1)g. \ea By (\ref{smallRmC-1}), we have now \ba
\label{smallintegralC-1} \int_0^T
\frac{1}{V(x_0,1,t)}\int_{B(x_0,1,t)}
|Rm|^{\frac{n+2}{2}}dqdt  \le \delta_0. \ea
 By Bishop-Gromov relative volume comparison, we have
\ba \label{volumeC-1}
V(x, R, t) \le \frac{V_{-1}(R)}{V_{-1}(r)} V(x, r, t)
\le C(n) \frac{V(x, r, t)}{r^n},
\ea
with a positive constant $C(n)$ depending only on $n$, provided that $t \in [0, T],
d(x_0,x, t)<1$, and $0<r<R\le 1-d(x_0,x,t)$. Here $V_{-1}(r)$ denotes the volume
of a geodesic ball of radius $r$ in $\mbox{\bf H}^n$, the $n$-dimensional
hyperbolic space
(of sectional curvature $-1$). If $t \in [0, T], d(x_0, x,t) \le \frac{1}{4}$,
we then have $B(x_0, \frac{1}{4}, t) \subset B(x, \frac{1}{2}, t)
\subset B(x_0, 1, t)$. Consequently,
\ba
V(x, r,t) &\ge& C(n)^{-1} r^n V(x, \frac{1}{2}, t)
\ge C(n)^{-1} r^n V(x_0, \frac{1}{4}, t) \nonumber \\  &\ge&
4^{-n} C(n)^{-2} r^n V(x_0, 1, t)
\ea
for $0<r\le \frac{1}{2}$.
Hence we infer
\ba \label{smallintegralC-1-2}
 \frac{r^n}{V(x,r,t)}\int_{B(x, r, t)} |Rm|^{\frac{n}{2}} dq \le
 \frac{4^n C(n)^2}{V(x_0, 1, t)}
 \int_{B(x_0, 1, t)}
|Rm|^{\frac{n}{2}} dq\ea whenever $0\le t <T, 0<r\le \frac{1}{2}$,
and $d(x_0, x, t) \le \frac{1}{4}$. This leads to 
\ba \label{Cvolume}
V(x, r, t) \ge 4^{-n}C(n)^{-2} \frac{r^n}{r_0^n}V(x_0, r_0, t)
\ea
as along as $t \in [0, T), d(x_0, x ,t) \le \frac{1}{4}r_0$ 
and $0<r\le \frac{1}{2}r_0$. 
Hence
 \ba
\label{smallintegralC-1-3} &&\int_0^T \sup\limits_{x \in B(x_0, \frac{r_0}{4}, t)} \frac{r^n}{V(x,r,t)}\int_{B(x, r, t)}
|Rm|^{\frac{n}{2}} dq dt \nonumber \\ &\le&  4^n C(n)^2
\int_0^T  \frac{r_0^n}{V(x_0, r_0, t)}
 \int_{B(x_0, r_0, t)}
|Rm|^{\frac{n}{2}} dq dt \le 4^n C(n)^2 \delta_0 \ea for all $x \in M$ 
and $0<r\le \frac{r_0}{2}$. 
Choosing $\delta_0$ to be the $\delta_0$ in Theorem B-1 multiplied
by $4^{-n}C(n)^{-2}$ and replacing $r_0$ by $\frac{r_0}{2}$ we
then have all the conditions of Theorem B-1. The desired estimates
follow (\ref{RmboundC-1}) and (\ref{RmboundintC-2}) follow. The 
estimate (\ref{RmboundintC-1}) follows from (\ref{RmboundintC-2}) and 
(\ref{Cvolume}). 
\\
\qed \\

\sect{Convergence}
\vspace{3mm}

\subsection{Type A}
\vspace{3mm}

\noindent {\bf Proof of Theorem A-2} \\

By [Theorem 4.1, P] or [Theorem A.1, Ye3], $g$ is $\kappa$-noncollapsed on the scale $\sqrt{T}$ for some $\kappa>0$ depending on 
$T$ and $g(0)$.  By (\ref{finiteRmA-2}) we can choose $0<T_0<T$ such that 
\ba \label{small-1}
\int_{T_0}^T\int_M |Rm|^{\frac{n+2}{2}}dqdt\le \delta_0,
\ea
where $\delta_0=\delta_0(\kappa, n)$ is from Theorem A-1. Then we have
\ba
\int_{T_0}^T\int_{B(x_0, \sqrt{T}, t)} |Rm|^{\frac{n+2}{2}} dqdt \le \delta_0
\ea
for all $x_0 \in M$. Obviously, $d_{g(t)}(x, \partial M)=\infty$ for all $x \in M$ and $ t\in [0, T)$. 
Hence we can apply Theorem A-1 with $r_0=\sqrt{T}$ and $T_0$ playing the role of the time origin $0$ to deduce
\ba \label{curvature-1}
|Rm|(x, t) \le \alpha_n 2(T-T_0)^{-1}+\epsilon_0^{-2} T^{-1}
\ea
for all $x \in M$ and $(T_0+T)/2\le t <T$. Since $T$ is finite, the desired smooth convergence follows. 
(Higher order estimates for $Rm$ follow from [Sh].  A local positive lower bound for volume follows 
from (\ref{curvature-1}) and the Ricci flow equation, or from (\ref{curvature-1}) and the $\kappa$-noncollapsedness. 
Then an injectivity radius estimate follows from [CGT] or [Lemma B.1, Ye3].) 
 \qed \\
 
 \noindent {\bf Proof of Theorem A-3} \\
 
 We obtain local curvature estimates and local injectivity radius estimates in the same way as in 
 the proof of Theorem A-2. The desired smooth convergence follows.
 The identification of the limit domain follows from an estimate of distance change
based on the Ricci flow equation and the obtained curvature estimate.  \\ \qed \\

\subsection{Type B} 

\vspace{3mm}

\noindent {\bf  Proof of Theorem B-2} \\

By [Theorem 4.1, P] or [Theorem A.1, Ye3], $g$ is 
$\kappa$-noncollapsed on the scale $\sqrt{T}$ for some $\kappa>0$ depending 
on $T$ and $g(0)$. 
By the evolution equation for the scalar curvature 
\ba
\frac{\partial R}{\partial t}=\Delta R+2|Ric|^2
\ea
and the maximum principle we have 
\ba \label{Rbound}
R_{min}(t) \ge R_{min}(0)
\ea
for all $0\le t<T$, where $R_{min}(t)$ denotes the minimum of $R$ at time $t$. 
Next observe that by (\ref{finiteRmB-2}) we can choose $T_0<T$ such that (with the $\delta_0=\delta_0(n, \kappa)$
from Theorem B-1)
\ba \label{smallB-2}
\int_{T_0}^T\sup_{x\in M, 0<r \le d(r_0, \mu g(t))}
\frac{r^n}{V(x, r, t)} \int_{B(x, r, t)}
|Rm|^{\frac{n+2}{2}}dqdt \le \delta_0.  \ea

 By rescaling, we can assume that the volume of $g(0)$
is 1. We rescale $g(t)$
to obtain a solution $\bar g(\tau)$ of the volume normalized Ricci
flow  on $M \times [0, \Lambda)$ with $\bar g(0)=g(0)$, where $\Lambda$
corresponds to $T$. Thus $\bar g(\tau)=\phi(t(\tau)) g(t(\tau))$, where
$\phi(t)=exp(\frac{2}{n}\int_0^t \hat R)$, $\tau(t)=\int_0^t \phi$, $t(\tau)$ is the inverse of 
$\tau(t)$, and 
$\hat R$ denotes the average of $R$.  Let $\Lambda_0$ correspond to $T_0$. i.e. $\Lambda_0=\tau(T_0)$.
By (\ref{Rbound}) 
we have
\ba \label{c-0}
\phi(t) \ge e^{\frac{2}{n}\int_0^t R_{min}} \ge  c_0 \equiv e^{\frac{2}{n} \min\{R_{min}(0), 0\}T}.
\ea
It follows that 
$\bar g$ is $\kappa$-noncollapsed on the scale $\sqrt{c_0T}$. 
By (\ref{smallB-2})
we infer for $\bar g$ (all quantities are associated with $\bar g$)  \ba \label{smallafter}
\int_{\Lambda_0}^{\Lambda}\sup_{x\in M, 0<r\le d( r(\tau),  \mu \bar
g(\tau)) } \frac{r^n}{V(x, r, \tau)} \int_{B(x, r,
\tau)} |Rm|^{\frac{n+2}{2}} dqd\tau\le \delta_0, \ea
where $r(\tau)=\phi^{\frac{1}{2}}(t(\tau)) r_0$.  On account of (\ref{c-0}) we have 
\ba \label{r-bound}
r(\tau) \ge r_1 \equiv \sqrt{c_0}r_0.
\ea

The scalar curvature $R_{\bar g}$ of $\bar g$ satisfies
\ba \label{barR}
R_{\bar g}=\phi^{-1}R \ge \phi^{-1}\min\{R_m(0), 0\} 
\ge c_0^{-1}\min\{R_m(0), 0\} =-a_0,
\ea
where $a_0=-c_0^{-1}\min\{R_m(0), 0\}$. 

We apply Theorem B-1 to derive a $|Rm|$ estimate for $\bar g$. For this purpose, consider
an arbitary interval $[\bar \tau-\sigma, \bar \tau]$ contained in $ [\Lambda_0, \Lambda)$
  We convert $\bar g$ on $[\bar \tau-\sigma, \bar \tau]$  into a solution $g^*(s)$ of the Ricci flow by rescaling. Namely we set
$g^*(s)=\psi(\tau(s)) \bar g(\tau(s))$, where $\psi(\tau)=\exp(-\frac{2}{n} \int_{\bar \tau}^{\tau} \int_M
\hat R_{\bar g})$, $s=\int_{\bar \tau-\sigma}^{\tau} \psi$, $\tau(s)$ is the inverse 
of $s(\tau)$, and $\hat R_{\bar g}$ denotes the average of $R_{\bar g}$.   We set $S^*=s(\bar \tau)$,
whence $g^*$ is defined on $[0, S^{*}]$. We infer from (\ref{smallafter}) the following 
estimate for $g^*$ (all quantities are associated with $g^*$)
\ba \label{smallonce}
\int_{\bar \tau-\sigma}^{\bar \tau}\sup_{x\in M, 0<r\le d( r^*(s), 
\mu g^*(s)) } \frac{r^n}{V(x, r, s)} \int_{B(x, r,
s)} |Rm|^{\frac{n+2}{2}}dqds \le \delta_0,
\ea
where $r^*(s)=\psi^{\frac{1}{2}}(\tau(s))r(\tau(s))$.

 Applying the estimate
(\ref{RmboundB-1}) in Theorem B-1 on the interval $[0, S^*]$
 we deduce for $g^*$ the
estimate 
\ba \label{gstar} |Rm(x, S^*)| \le \alpha_n
({S^*})^{-1}+(\epsilon_0 d^*)^{-2}, \ea where \ba \label{d*} d^*=\min_{[0,
S^*]} d(r^*(s), \mu g^*(s))= \min_{[\bar \tau-\sigma, \bar \tau]} \psi^{\frac{1}{2}}(\tau) 
d(r(\tau), \mu \bar g(\tau))
\ea 
for all $x\in M$.
Scaling back we obtain the following
estimate for $\bar g$ at time $\bar \tau$ 
\ba \label{barR-1} 
|Rm|(x, \bar
\tau) \le \alpha_n \frac{\psi(\bar \tau)}{S^*} + \psi(\bar
\tau)(\epsilon_0 d^*)^{-2} \ea
for all $x \in M$.

We have 
\ba \label{psi}
\psi(\bar \tau) \le e^{\frac{2}{n}a_0 \sigma} \psi(\tau)
\ea
for all $\tau \in [\bar \tau-\sigma, \bar \tau]$. 
 Since $S^*=\int_{\bar
\tau-\sigma}^{\bar \tau} \psi$, we infer
\ba \label{star} 
S^* \ge \sigma e^{-\frac{2}{n}a_0 \sigma} \psi(\bar \tau). 
\ea
By (\ref{psi}) and (\ref{d*})  we also deduce 
\ba \label{bar-d}
 d^* 
 \ge e^{-\frac{1}{n}a_0 \sigma}\psi^{\frac{1}{2}}(\bar \tau) \bar d, 
 \ea
 where 
 \ba 
 \bar d =\min_{[\bar \tau-\sigma, \bar
\tau]} d(r(\tau), \mu \bar g(\tau)). \ea 

By (\ref{barR-1}), (\ref{star}) and (\ref{bar-d}) we arrive at the
following estimate for $\bar g$
 \ba \label{analog1}
|Rm|(x, \bar \tau) \le e^{\frac{2}{n}a_0 \sigma}(\alpha_n \sigma^{-1}+
(\epsilon_0 \bar d)^{-2}) \ea
for all $x \in M$. \\

\noindent {\bf Claim} $\liminf\limits_{\tau \rightarrow \Lambda} diam_{\bar g(\tau)}(M)>0.$\\

Assume the contrary. Then we can find a sequence $\tau_k \rightarrow \Lambda, \tau_k
> \Lambda_0$, such that $d_k\equiv diam_{\bar g(\tau_k)}(M) \rightarrow 0$
and $d_k=\inf\limits_{\Lambda_0\le \tau \le \tau_k}
diam_{\bar g(\tau)}(M)$.  Now we apply (\ref{analog1}) to the interval
$[\bar \tau-\sigma, \bar \tau]$ with $\bar  \tau=\tau_k$ and $\sigma=d_k^2$ (for $k$ large enough) to deduce
\ba
|Rm|(x, \tau_k) \le e^{\frac{2}{n}a_0 d_k^2}(\alpha_n d_k^{-2}+
(\epsilon_0 \bar d)^{-2}).
\ea
By (\ref{r-bound}) and the properties of $d_k$ we  have $\bar d=\sqrt{\mu} d_k. $  Hence we infer
\ba
|Rm|(x, \tau_k) \le e^{\frac{2}{n}a_0 d_k^2}(\alpha_n+
\epsilon_0^{-2} \mu^{-1}) d_k^{-2}.
\ea

 Now the rescaled metric $d_k^{-2} \bar g(\tau_k)$ has
diameter 1 and satisfies $|Rm| \le e^{\frac{2}{n}a_0 d_k^2}(\alpha_n+
\epsilon_0^2 {\mu})$. Since $d_k \rightarrow 0$, we obtain 
$|Rm| \le \bar C$ for a positive constant $\bar C$ independent of $k$.  By volume comparison,  the volume of 
$\bar g(\tau_k)$ is bounded
from above by a constant independent of $k$. On the other hand, its volume equals $d_k^{-\frac{n}{2}}
$,
which approaches $\infty$ as $k \rightarrow \infty$. This is a contradiction. Thus the claim is proved. \\

By (\ref{analog1}) and the above Claim we deduce for $\bar g$
\ba \label{analog2}
|Rm|(x, \tau) \le C
\ea
for a positive constant $C$, all $x\in M$ and all $\tau \in [\Lambda_0, \Lambda)$. \\

Next we employ the estimate 
(\ref{RmboundintB-1}) in Theorem B-1 to derive
the following curvature estimate for the above $g^*$ 
(all quantities are associated with $g^*$)
\ba 
|Rm|(x, S^*) \le C_0 r_*^{-2}
\left(\int_0^{S^*} \frac{r_*^n}{V(x, r_*, s)} \int\limits_{B(x, r_*, s)} |Rm|^{\frac{n+2}{2}}dqds
\right)^{\frac{2}{n+2}}
\ea  
where $r_*=\frac{1}{9}\epsilon_1 \min\{d^*, \sqrt{S^*}\}$. 
We convert this estimate into an estimate for $\bar g$  (all quantities are associated with $\bar g$)
 \ba \label{int-1} && |Rm|(x, \bar \tau) \psi(\bar \tau + \sigma)^{-1} \le  \nonumber \\
 && C_0 \psi(\bar \tau) r_*^{-2}
\left(\int_{\bar \tau-\sigma}^{\bar \tau} \frac{(\psi^{-\frac{1}{2}}(\tau)r_*)^n}{V(x, \psi^{-\frac{1}{2}}(\tau)r_*, \tau)}\int\limits_{B(x, \psi^{-\frac{1}{2}}(\tau)r_*,
\tau)} |Rm|^{\frac{n+2}{2}}dqd\tau
\right)^{\frac{2}{n+2}}.
\ea 
By (\ref{star}) and (\ref{bar-d}) we infer
\ba 
\psi(\bar \tau) r_*^{-2} \le 81\epsilon_1^{-2} e^{\frac{2}{n}a_0\sigma} \max\{\bar d^{-2}, 
{\sigma}^{-1}\}.
\ea
On the other hand, we have by (\ref{d*}) 
\ba
\psi^{-\frac{1}{2}}(\tau) r_* \le \frac{1}{9} \epsilon_1 \psi^{-\frac{1}{2}}(\tau) d^* \le 
\frac{1}{9} \epsilon_1 d(r(\tau), \mu \bar g(\tau)) \le  d(r(\tau), \mu \bar g(\tau)).
\ea

Hence we infer for $\bar g$  
\ba \label{analog3} &&|Rm|(x, \bar \tau)  \le   81\epsilon_1^{-2}e^{\frac{2}{n}a_0\sigma} \max\{\bar d^{-2}, \sigma^{-1}\} \nonumber \\ 
&& \cdot 
\left(\int_{\bar \tau}^{\bar \tau+\sigma} \sup\limits_{0<r\le d(r(\tau), \mu \bar g(\tau))} 
\frac{r^n}{V(x,r,\tau)}\int\limits_{B(x,r,\tau)} |Rm|^{\frac{n+2}{2}}dq d\tau
\right)^{\frac{2}{n+2}}. 
\ea

Now we divide into two possible cases. \\

\noindent {\bf Case 1} $\Lambda<\infty$. \\

We convert $\bar g$ on $[\Lambda_0, \Lambda)$ back into $g$ in the same 
way as converting $\bar g$ on $[\bar \tau-\sigma, \bar \tau]$ into $g^*$.
Namely we have $g(t)=f(\tau(t)) \bar g(\tau(t))$, where 
$f(\tau)=\exp(-\frac{2}{n}\int_{\Lambda_0}^{\tau} \hat R_{\bar g})$,
$t=\int_{\Lambda_0}^{\tau} f$, and  $\tau(t)$ denotes the inverse of 
$t(\tau)$.  Since $\Lambda<\infty$, we deduce from (\ref{analog2}) 
that $f\ge c$ for a positive constant $c$ independent of $\tau$.  Since the $|Rm|$ of $g$ is given by 
the $|Rm|$ of $\bar g$ multiplied by $f^{-1}$, we obtain a uniform
$|Rm|$ bound for $g$ over $[0, T)$. It follows that
$g(t)$ converges smoothly as $t \rightarrow T$. \\

\noindent {\bf Case 2} $\Lambda=\infty$. \\

By (\ref{smallafter}) we have 
\ba \label{smallatinfty}
\lim\limits_{\tau^* \rightarrow \infty} \int_{\tau^*}^{\infty}\sup_{x\in M, 0<r\le d( r(\tau),  \mu \bar
g(\tau)) } \frac{r^n}{V(x, r, \tau)} \int_{B(x, r,
\tau)} |Rm|^{\frac{n+2}{2}}dqd\tau=0.
 \ea

We apply (\ref{analog3}) for $\bar \tau \ge \Lambda_0+1$ and $\sigma=1$. On account of 
(\ref{smallatinfty}) and the above Claim
we then deduce 
\ba \label{zerocurvature}
\limsup_{\tau \rightarrow \infty} \sup\limits_{x\in M} |Rm|(x,  \tau)=0.
\ea
Since $\bar g$ is $\kappa$-noncollapsed on the scale $\sqrt{c_0T}$, we conclude that $\bar g(\tau)$ subconverges smoothly to flat metrics
 of volume 1 on $M$ as $\tau \rightarrow \infty$. For each
$\bar t>0$ we rescale $g$ by the constant factor $\phi(\bar t)$ to obtain
 $g^{\bar t}(t)=
\phi(\bar t) g(\phi(\bar t)^{-1}t)$ on $[\phi(\bar t)\bar t, \phi(\bar t)T)$.
Note that $g^{\bar t}$ is a solution of the Ricci flow with $g^{\bar t}(\phi(\bar t) \bar t)=\bar
g(\tau(\bar t))$. As $\bar t \rightarrow T$, we have $\tau(\bar t) \rightarrow \infty$,
and hence $g^{\bar t}(\phi(\bar t) \bar t)$ subconverges smoothly to flat metrics of  volume 1
on $M$.
By the stability theorem in [GIK], each flat metric has a smooth neighborhood such that
the Ricci flow starting at a metric in the neighborhood exists for all time and converges
smoothly to a flat metric at the time infinity.
It follows that $g^{\bar t}(t)$ extends to a smooth solution of the Ricci flow
 for all time $t \in  [\lambda(\bar t) \bar t, \infty)$ and converges smoothly to a flat metric
as $t \rightarrow \infty$, provided that $\bar t$ is close enough to
$T$.  Consequently, $g(t)$ converges smoothly as $t \rightarrow T$, and
$g$ extends to $[0, \infty)$ and converges to a flat metric as $t \rightarrow \infty$.
\\

\noindent {\bf Proof of Theorem B-3} \\

This is similar to the proof of Theorem A-3. \qed

\subsection{Type C} 
\vspace{3mm}

\noindent {\bf Proof of Theorem C-2} \\

By the arguments in the proof of Theorem C-1, we can reduce to the
situation of Theorem B-2. \qed \\

\noindent {\bf Proof of Theorem C-3} \\

This is similar to the proof of Theorem A-3. \\
\qed

\end{document}